\documentclass[twoside]{article}
\def\Vox{$^{\fbox {\/}}$}
\def\prom{
\newtheorem{prop}{Proposition}[section]
\newtheorem{cor}[prop]{Corollary}
\newtheorem{lem}[prop]{Lemma}

\newtheorem{rem}[prop]{Remark}

\newtheorem{tef}[prop]{Definition}
}
\def\ieacosdo{
\author {\ \\N.\ C.\ A.\ da Costa\\F.\ A.\ Doria\\\ 
\\Institute for Advanced Studies,
University of S\~ao Paulo.\\Av.\ Prof.\ Luciano
Gualberto,
trav.\ J, 374.\\05655--010 S\~ao Paulo SP Brazil.\\\
\\{\sc
ncacosta@usp.br}\\{\sc doria1000@yahoo.com.br}\\{\sc
doria@lncc.br}}}
\prom

\begin {document}
\pagestyle {myheadings}
\title {A lemma on a total function defined over the
Baker--Gill--Solovay set of polynomial Turing
machines.\thanks {Partially supported by FAPESP and
CNPq.
Alternative address for F.\ A.\ Doria: Research Center
for
Mathematical Theories of Communication and Program
IDEA,
School of Communications, Federal University at Rio de
Janeiro, Av.\ Pasteur, 250. 22295--900 Rio RJ
Brazil.}}
\ieacosdo
\date {\today}
\maketitle

\begin {abstract}
\noindent We prove here a lemma that connects some
properties of the so--called ``counterexample
function'' to
the $P=NP$ conjecture over the Baker--Gill--Solovay
set of
polynomial Turing machines to the behavior of the same
function ``at large,'' over the set of all polynomial
Turing
machines. 
\end {abstract}

\newpage
\thispagestyle {plain}

\markboth {da Costa, Doria}{BGS set}
\section {Introduction}
\markboth {da Costa, Doria}{BGS set}

We deal here with a property of the
Baker--Gill--Solovay
(BGS) set \cite {Bak} of polynomial Turing machines.
The BGS
set is a kind of `representation set' for polynomial
machines in the following sense: for every computable
function ${\sf f}$ with a polynomial algorithm (a
polynomial
Turing machine), there is one such polynomial
algorithm for
${\sf f}$ in the BGS set, and there are only
polynomial
algorithms in BGS. 

Actually there are infinitely many polynomial
algorithms for
each ${\sf f}$ in BGS, but not every polynomial
algorithm for
${\sf f}$ will be in BGS. 

The BGS set was conceived to rigorously formulate the
$P<NP$
question (see below). Since the set of all polynomial
Turing
machines isn't recursive in the set of all Turing
machines,
it is certainly easier to deal with a set that
contains
copies of representatives of all polynomial algorithms
and
which moreover is a recursive set. 

\subsection*{On our main result}

As we see below, $P<NP$ asserts that a given recursive
function (noted ${\sf f}_{\neg G}$ in what follows) is
total
over the BGS set. That function has an extension
(better
said, a kind of copy) as a relative recursive function
on
the set of all polynomial Turing machines. So, it is
of
interest to relate what happens over the BGS set to
what
happens outside it. 

We are interested in Peano Arithmetic (PA), which we
take to
be consistent. We deal with PA--provably total
recursive
functions. Recall that:

\begin {tef}
${\sf F}$ is {\rm PA}--provably total recursive if, 
\begin {enumerate}
\item ${\sf F}$ has an explicitly given G\"odel number
$e$. 
\item ${\rm PA}\vdash\forall x\,\exists z\,T(e,x,z)$.
\Vox
\end {enumerate}
\end {tef}

($T$ is Kleene's predicate.) This means that we
explicitly
have a program for ${\sf F}$ (it is given by $e$) and
that
there is a proof in PA that every computation of ${\sf
F}$ 
will eventually stop. 

Our question in this paper is: suppose that the
counterexample function is total over BGS. If it is
PA--provably total, what happens outside BGS, over the
(nonrecursive) set of all polynomial Turing machines?

We give here a partial answer to that question. 

\subsection* {The BGS set}

The BGS set is constructed as follows: 
\begin {itemize} 

\item A polynomial clock ${\sf C}_{(a,b)}$ is a total
Turing machine that behaves as follows: for binary
input $x$
of length $|x|$ it computes $|x|^a + b$, $a, b$
positive
integers, and stops the operation of the coupled
machine
${\sf M}_m(x)$  after $|x|^a + b$ cycles, if it hasn't
stopped yet. 

\item Form all pairs $({\sf M}_m, {\sf C}_{(a,b)})$ of
a
Turing machine ${\sf M}_m$ and a polynomial clock
${\sf
C}_{(a,b)}$. 

\item The pairs $({\sf M}_m, {\sf C}_{(a,b)})$ form
the BGS
set. If $\langle\ldots,\ldots\rangle$ is the usual
1--1 and
onto pairing function, we order BGS according to
$\langle m,
\langle a, b\rangle\rangle$. 

The BGS index is the triple $\langle m,a,b\rangle$,
$m, a, b$ ranging over the whole of $\omega$.

\item There are several primitive recursive procedures
to 
embed all such pairs into the set ${\cal M}$ of all
Turing
machines. We suppose that one of them has been chosen
and
kept fixed.  The pairs $({\sf M}_m, {\sf C}_{(a,b)})$
or
their recursive, embedded images, form the BGS set.
(For our
purposes it is indifferent whether we deal with the
BGS
pairs or with their image in ${\cal M}$, the set of
all
Turing machines, via that p.r.\ embedding previously
agreed
upon, but we will consider here the BGS set as an
entity
separated from ${\cal M}$.) 

\end {itemize}

As mentioned above, one immediately sees that for
every poly
machine there is a pair machine--clock in BGS that
corresponds to it (actually, infinitely many such
pairs),
and given an arbitrary pair, there is a corresponding
poly
machine. \Vox

\subsection*{Rigorous formulation of $P<NP$ for the
Satisfiability Problem}

\begin {rem}\label {AAA}\rm We consider the case \cite
{CosDo0} of {\sc Sat}, the  Satisfiability Problem for
Boolean
expressions in conjunctive normal form (cnf). 
\begin {itemize}

\item Let $x$ be a Boolean expression in cnf,
adequately
coded as a binary string of length
$|x|$. Let ${\sf P}_n$ be a polynomial machine of BGS
index
$n = \langle m, a, b\rangle$. 

\item Given a binary string $y$ of truth--values for
the
$|y|$ Boolean variables of $x$, there is a polynomial
procedure (a polynomial Turing machine which we note
${\sf V}$) that tests whether $y$ satisfies $x$, that
is,
say, ${\sf V}(\langle x,y\rangle)=1$ if and only if
$y$
satisfies $x$; and is $0$ otherwise. 

(For the sake of completeness, we add that ${\sf
V}(0,0)=1$,
that is, the empty string is satisfied by the empty
string. 
The empty string as a string of truth values makes
true the
empty string, seen as a string of propositional
variables;
we agree that for $x>0$, no such $x$ is satisfied by
$0$.)  

\item We formulate the predicate: 
$$G^*(m,x)\leftrightarrow _{\rm Def}\exists y\,({\sf
P}_m(x)=y\,\wedge\,{\sf V}(x,y)=1).$$ 
$G^*(m,x)$ is intuitively understood as ``polynomial
machine
of BGS index $m$ correctly guesses about Boolean cnf
expression $x$,'' or, even more explicitly, ``machine
$m$
inputs $x$ and outputs a line of truth values that
satisfies
$x$.''

\item We can also write: 
$G^*(m,x)\leftrightarrow\,[{\sf V}(x,{\sf
P}_m(x))=1]$. 

\item Form the pair $z=\langle x,y\rangle$, and let
$\pi_i$,
$i=1,2$, be the usual (polynomial) projection
functions.
Recall that ${\sf V}$ is a polynomial machine that
inputs a
pair
$\langle x,y\rangle$. Then we can consider the
predicate: 
$$\neg G(m, z)\leftrightarrow _{\rm Def}{\sf
V}(z)=1\,\wedge\, {\sf V}(\langle \pi _1 z, {\sf P}_m
(\pi
_1 z)\rangle)=0,$$  or
$$\neg G(m,z)\leftrightarrow {\sf V}(z) =
1\,\wedge\,\neg
G^*(m,\pi_1 z).$$ 

\item $\neg G(m,z)$ can be intuitively understood as
follows:
polynomial machine ${\sf P}_m$ doesn't accept the pair
$z$
if and only if $z$ is such that $\pi _1 z = x$ is
satisfiable, but the output of ${\sf P}_m$ over $x =
\pi_1
z$ doesn't satisfy $x$. \Vox
\end {itemize}
\end {rem} 

Then we can define:

\begin {tef}\label {p<np}
$P < NP\,\leftrightarrow_{\rm Def} \forall m\,\exists
z\,\neg G(m, z)$, where $m$ ranges over the BGS set.
\Vox
\end {tef}

Notice that $P<NP$ is a $\Pi^0_2$ sentence. Also:

\begin {tef} ${\sf f}_{\neg G} =_{\rm Def} \mu_x [\neg
G(m,
x)]$ is the {\bf counterexample function}. \Vox
\end {tef}

\markboth {da Costa, Doria}{BGS set}
\section{Main result}
\markboth {da Costa, Doria}{BGS set}

\begin {rem}\rm  We suppose here that the
counterexample
function ${\sf f}_{\neg G}$ is total over BGS. This
means
that we suppose that $P<NP$ holds, for the sake of our
argument.

Our query then is: if it is so, what do we need in
order to
have that ${\sf f}_{\neg G}$ be PA--provably total
recursive? \Vox
\end {rem}

Write $f^*_{\neg G}$ for the `complete' counterexample
function, that is the one which is defined over the
(nonrecursive) set ${\cal P}\subset {\cal M}$ of all
poly
Turing machines. 

Let ${\sf T}$ be the exponential algorithm
(truth--table
computation) that settles {\sc Sat}. We have agreed
that no
instance $> 0$ is satisfied by $0$. 

\begin {rem}\rm 
An ${\sf A}$--{\it
quasi--trivial Turing machine with cutoff value $k$}
is a
Turing machine that equals some other total machine
${\sf A}$
up to instance $k$, and then outputs $0$ for every
instance
$x>k$.
\Vox
\end {rem}

Define:

\begin {tef} A ${\sf T}${\bf --quasi--trivial machine}
${\sf
T}^k$ is a quasi--trivial machine that equals ${\sf
T}$ up
to instance $k$, and then outputs $0$ for every
instance
$x>0$. \Vox
\end {tef}

\begin {tef} A recursive subset $B$ contained in the
set of
all Turing machines is {\rm PA}--{\bf recursive} iff
its
explicitly given characteristic function ${\sf c}_B$
is {\rm
PA}-provably total recursive. \Vox
\end {tef}

\begin {rem}\rm
If $B$ as above is PA--provably recursive, then we can
also
explicitly obtain a PA--provably recursive function
${\sf
c}_B: B\rightarrow\omega$ that is 1--1 and onto. Thus
we can
code (via ${\sf c}_B$) the machines in $B$ by the
natural
numbers. 

Given a PA--provably total recursive function ${\sf
F}$ which is defined over all Turing machines coded by
$\omega$, we can adequately define the restriction
${\sf
F}|_B$ for $B$ PA--provably total recursive, and see
that ${\sf
F}|_B$ is PA--provably total over the set $B$ coded by
${\sf
c}_B$. \Vox
\end {rem}

\begin {rem}\rm \label {N}
Recall that, for $m$ the G\"odel number of a
polynomial
machine, and for the corresponding BGS index $N(m)$,
$$f^*_{\neg G}(m) = {\sf f}_{\neg G}(N(m)).$$ 

The map $m\mapsto N(m)$ isn't in general recursive (it
isn't
even a function in the general case, as there are
infinitely many $N(m)$ that correspond to each $m$).
However
we use below recursive versions of it which also turn
out to
be functions. \Vox

\end {rem}

We now state our main result: 

\begin {lem} If the counterexample function ${\sf
f}_{\neg G}$ is {\rm PA}--provably total over BGS,
then for any
restriction of $f^*_{\neg G}$ over a {\rm
PA}--recursive
subset $B$ of {\sf T}--quasi--trivial machines,
$f^*_{\neg
G}|B$ is {\rm PA}--provably total recursive. \Vox
\end {lem} 

{\it Proof of the lemma}\,: Keep in mind the
correspondence $$f^*_{\neg G}(m) = {\sf f}_{\neg
G}(N(m))$$ given in Remark \ref {N}. Pick up an
arbitrary
PA--recursive subset $B$ of the Turing machines which
only
contains ${\sf T}$--quasi--trivial machines. Let's
embed it
into BGS as follows:
\begin {itemize}

\item {\it Clocks that bound $B$.} Machines in $B$ are
as
follows ($m$ is the machine's G\"odel number):
$${\sf T}_m^{k(m)}(x) = {\sf T}(x), x\leq k(m),$$
$${\sf T}_m^{k(m)}(x) = 0, x>k(m).$$

\item Put 
$b_m = \max\{\mbox{operation time of ${\sf T}(x),
x\leq
k(m)$}\} + 1$. 

\item Then clock ${\sf C}_{(2,b_m)}$ bounds the
operation of
${\sf T}_m^{k(m)}$ without interrupting it. 

\item We thus form $B'\subset {\rm BGS}$, whose
elements are
the pairs $({\sf T}_m^{k(m)}, {\sf C}_{(2,b_m)})$.
This
gives us the recursive map (see Remark \ref {N}):
$$B\subset {\cal P}\rightarrow B'\subset {\rm BGS},$$
$$m\mapsto N(m) = \langle m, 2, b_m\rangle.$$

\item {\it Crucial step.} Now we know that ${\sf
f}_{\neg
G}$ is PA--provably recursive over BGS. For each pair
$({\sf
T}_m^{k(m)}, {\sf C}_{(2,b_m)})$, and if the BGS index
$N(m)
= \langle m, 2, b_m\rangle$, then:
$${\sf f}_{\neg G}(N(m))\geq k(m) + 1,$$ by the
definition of ${\sf f}_{\neg G}$. 

\item $m\in B$. Since $B$ is PA--provably recursive,
this
means that the $k(m)$ are also PA--provably recursive
over
$B'$. 

\item Therefore, so are the $b_m$ by construction, and
as a
result $B'$ has a PA--provably recursive
characteristic
function in BGS. Thus $B'$ is PA--recursive as a
subset of
BGS. 

\item {\it Conclusion.} Now go back to $B$ and
trivially
obtain the values of $f^*_{\neg G}|B$ from those of
${\sf
f}_{\neg G}|B'$. 

Proceed as follows: for $\langle m, 2,
b_m\rangle \in B'$ and $m\in B$, we have (Remark
\ref {N}) that: 
$$f^*_{\neg G}(m)|_B = {\sf f}_{\neg G}(\langle m, 2,
b_m\rangle). \mbox {\Vox}$$ 
\end {itemize}

\begin {cor} If there is a restriction $f^*_{\neg
G}|B$
which isn't PA--provably total, then ${\sf f}_{\neg
G}$
(over BGS) cannot be PA--provably total. \Vox
\end {cor}

Therefore, if we manage to show that at least one such
restriction isn't PA--provably total recursive, we
have that
the counterexample function over BGS cannot be
PA--provably
total recursive. Consequence is:

\begin {cor}
If there is a restriction $f^*_{\neg G}|B$ which isn't
PA--provably total, then PA cannot prove $P<NP$. 
\end {cor}

{\it Proof}\,: Follows from the fact that $P<NP$ is a
$\Pi^0_2$ sentence in PA, and from Kreisel's theorem
(\cite
{Schwich}, p.\ 885ff). \Vox

\

On that last possibility see \cite {CosDo1,CosDo0,te}.

\markboth {da Costa, Doria}{BGS set}
\section {Acknowledgments}
\markboth {da Costa, Doria}{BGS set}

This note is part of the Research Project on
Complexity and
the Foundations of Computation at the Institute for
Advanced
Studies, University of S\~ao Paulo (IEA--USP). 

Most of the discussion that led to this note was
carried out
in the forum {\bf theory--edge} at Yahoo Groups from
November to December 2000; please search its posts and
files
for more detailed references and due acknowledgments.
We are
very grateful for the comments, both public and
private, of
several participants of {\bf theory--edge}. 

The authors thank the Institute for
Advanced Studies, University of S\~ao Paulo, for
support, especially its Director Prof.\ A.\ Bosi, as
well as F.\ Katumi and S.\ Sedini. The second author
also
wishes to thank Professors S.\  Amoedo de Barros and
A.\
Cintra at the Federal University of Rio de Janeiro, as
well
its Rector J.\ Vilhena. 

\bibliographystyle {plain}

\end {document}